\documentstyle{article}
\newtheorem{prop}{Proposition}[section]
\newtheorem{th}[prop]{Theorem}
\newtheorem{cor}[prop]{Corollary}
\newtheorem{lm}[prop]{Lemma}

\newtheorem{remark}[prop]{Remark}

\newcommand{\Ker}{{\rm{Ker}}}
\newcommand{\bsquare}{\hbox{\rule{6pt}{6pt}}}
\newcommand{\proof}[1]{\noindent{\bf Proof}\hspace{0.3cm}{#1}\hfill\bsquare 
\vspace{0.5cm}\par}

\date{}

\begin{document}

\title{The Chern class map on abelian surfaces}
\author{Toshiyuki Katsura\thanks{Partially supported by JSPS Grant-in-Aid 
for Scientific Research (C) No. 24540053.}}

\maketitle

\begin{abstract}
We examine the Chern class map ${c}_{1}: {\rm NS}(S)/p{\rm NS}(S) \rightarrow {\rm H}^{1}(S, \Omega^{1}_{S})$ for an abelian surface $S$ 
in characteristic $p \geq 3$, and give a basis of the kernel $c_{1}$
for the superspecial abelian surface. 
\end{abstract}

\section{Introduction}
Let $k$ be an algebraically closed field of characteristic $p> 0$, and $S$ be 
a nonsingular complete algebraic surface over $k$. 
We denote by ${\rm H}^{2}_{dR}(S)$ the second de Rham cohomology group
of $S$, and by ${\rm NS}(S)$ the N\'eron-Severi group of $S$. ${\rm NS}(S)$
is a finitely generated abelian group, and the rank $\rho (S)$ of ${\rm NS}(S)$
is called the Picard number. 
We have the Chern class map ${\rm NS}(S)/p{\rm NS}(S) \rightarrow {\rm H}^{2}_{dR}(S)$
and this map is injective if the Hodge-to-de Rham spectral sequence of $S$
degenerates at $E_{1}$-term (cf. Ogus \cite{O}).
We also have the Chern class map 
${c}_{1}: {\rm NS}(S)/p{\rm NS}(S) \rightarrow {\rm H}^{1}(S, \Omega^{1}_{S})$.
This map is not necessarily injective, even if the Hodge-to-de Rham spectral 
sequence of $S$
degenerates at $E_{1}$-term (cf. Ogus \cite{O}).

In this paper, we examine this map $c_{1}$ in the case of abelian surfaces.
For abelian surfaces, the Chern class map $c_{1}$ is injective 
if and only if the abelian surface is not superspecial
(for the definition, see Section \ref{2}).  This fact was implicitly proved in Ogus \cite{O}
by using the notion of K3 crystal. We give here a down-to-earth proof of this fact
and  determine a basis of the kernel of the Chern class map $c_1$ 
for the superspecial abelian surface.
To calculate a basis of $\Ker~ c_{1}$, in Section \ref{2} we 
 examine the structure of the N\'eron-Severi group of
the superspecial abelian surface. Using theory of quaternion algebra,
problems on divisors on superspecial abelian surfaces are translated into
problems in matrix algebras over quaternion algebras.
As an example, we give an explicit description 
of our theory in the case of characteristic 3.
Finally, we examine the Chern class map for Kummer surfaces
and show results similar to those for abelian surfaces.

The author would like to thank Professor Gerard van der Geer
for useful advice and stimulating conversation. 
He is also grateful
to the referee for his/her careful reading and many suggestions.

\section{The N\'eron-Severi group}\label{2}
Let $k$ be an algebraically closed field of characteristic $p > 0$.
An abelian surface is said to be supersingular if it is isogenous to
a product of two supersingular elliptic curves.
An abelian surface is said to be superspecial if it is isomorphic to
a product of two supersingular elliptic curves. By definition,
if an abelian surface is superspecial, then the abelian surface is supersingular.
But the converse does not necessarily hold (cf. Oort \cite{Oo2}).
Note that
a superspecial abelian surface is unique up to isomorphism (cf. Shioda \cite{S}). 
In this section, we examine the structure of the N\'eron-Severi 
group of the superspecial abelian surface. 

Let $E$ be a supersingular elliptic curve defined over $k$, and
we consider the superspecial abelian surface $A = E_{1} \times E_{2}$
with $E_{1} = E_{2} = E$. We  denote by $O_{E}$
the zero point of $E$. We take a divisor $X = E_{1}\times\{O_{E_2}\} + 
\{O_{E_1}\} \times E_{2}$, which gives a principal polarization on $A$.
We also denote $E_{1}\times\{O_{E_2}\}$  
(resp.  $\{O_{E_1}\} \times E_{2}$) by $E_{1}$ (resp. by $E_{2}$)
for the sake of simplicity.
We set ${\cal O} = {\rm End}(E)$ and $B = {\rm End}^{0}(E)=
{\rm End}(E)\otimes {\bf Q}$.
Then $B$ is a quaternion division algebra over the rational number field ${\bf Q}$
with discriminant $p$, and ${\cal O}$ is a maximal order of $B$ (cf. Mumford \cite{M}, Section 22 and Deuring \cite{D}, Section 2).
For an element $a \in B$, we denote by $\bar{a}$ the image under the canonical involution.

We have a natural identification of ${\rm End}(A)$ with the ring ${\rm M_{2}}({\cal O})$
of two-by-two matrices with coefficients in ${\cal O}$:
$$
               {\rm End}(A) = {\rm M_{2}}({\cal O}).
$$
Here, the action of $\left(
\begin{array}{cc}
\alpha  & \beta \\
\gamma & \delta
\end{array}
\right) \in {\rm M_{2}}({\cal O})$ is given by
$$
\begin{array}{rccc}
\left(
\begin{array}{cc}
\alpha  & \beta \\
\gamma & \delta
\end{array}
\right) : &  A = E\times E & \longrightarrow & A = E\times E \\
    & (x, y)   &\mapsto & (\alpha (x) + \beta (y), \gamma (x) + \delta (y)).
\end{array}
$$

From here on, by a divisor $L$ we often mean the divisor class
represented by $L$ in ${\rm NS}(A)$ if confusion is unlikely to occur.
For a divisor $L$, we have a homomorphism
$$
\begin{array}{cccc}
   \varphi_{L} :  &A  &\longrightarrow  & {\rm Pic}^{0}(A) \\
      & x & \mapsto & T_{x}^{*}L - L,
\end{array}
$$
where $T_{x}$ is the translation by $x \in A$ (cf Mumford \cite{M}).
We set 
$$
     H = \{
\left(
\begin{array}{cc}
\alpha  & \beta \\
\gamma & \delta
\end{array}
\right) \in {\rm M}_{2}({\cal O})
~\mid
~\alpha, \delta \in {\bf Z},~\gamma, \beta \in {\cal O},~\gamma = \bar{\beta}
\}.
$$
The main part of the following theorem may be known to specialists
 (cf. Mumford \cite{M}, and Ibukiyama, Katsura and Oort \cite{IKO}), 
but since we cannot find a convenient reference, we give here a proof for it.
\begin{th}\label{intersection}
The homomorphism
$$
\begin{array}{cccc}
j : & {\rm NS}(A) & \longrightarrow  & H \\
   & L  &\mapsto  & \varphi_{X}^{-1}\circ\varphi_{L}
\end{array}
$$
is bijective.  By this correspondence, we have
$$
j(E_{1}) = \left(
\begin{array}{cc}
0  & 0 \\
0 & 1
\end{array}
\right),~
j(E_{2}) = \left(
\begin{array}{cc}
1  & 0 \\
0 &  0
\end{array}
\right).
$$
For $L_{1}, L_{2} \in {\rm NS}(A)$ such that
$$
  j(L_{1}) =
\left(
\begin{array}{cc}
\alpha_{1}  & \beta_{1} \\
\gamma_{1} & \delta_{1}
\end{array}
\right),~
  j(L_{2}) =
\left(
\begin{array}{cc}
\alpha_{2}  & \beta_{2} \\
\gamma_{2} & \delta_{2}
\end{array}
\right),
$$
the intersection number  $L_{1}\cdot L_{2}$ is given by
$$
L_{1}\cdot L_{2} = \alpha_{2}\delta_{1} + \alpha_{1}\delta_{2} - \gamma_{1}\beta_{2}
-\gamma_{2}\beta_{1}.
$$
In particular, for $L\in {\rm NS}(A)$ such that 
$j(L) = \left(
\begin{array}{cc}
\alpha  & \beta \\
\gamma & \delta
\end{array}
\right)
$
we have
$$
\begin{array}{l}
L^{2} = 2\det \left(
\begin{array}{cc}
\alpha  & \beta \\
\gamma & \delta
\end{array}
\right)\\
L\cdot E_{1}  = \alpha,~L\cdot E_{2} = \delta .
\end{array}
$$
We have also  $j(nD)=nj(D)$ for an integer $n$.
\end{th}
The first  and the final statements of this theorem are given 
in Mumford \cite{M}. In particular, the final statement follows easily 
from the definition of $\varphi_{L}$. 
To prove the others, we need some lemmas.
\begin{lm}
The restriction homomorphism
$$
\begin{array}{ccc}
{\rm Res} : {\rm Pic}^{0}(A) & \longrightarrow & {\rm Pic}^{0}(E_{1}) \times {\rm Pic}^{0}(E_{2}) \\
   L & \mapsto & (L\mid_{E_{1}}, L\mid_{E_{2}})
\end{array}
$$
is an isomorphism, and the following diagram commutes:
$$
\begin{array}{cccl}
A  & \stackrel{\varphi_{X}}{\longrightarrow} &  {\rm Pic}^{0}(A) & \ni \quad L\\
\mid\mid &   & \downarrow {\rm Res}  & \qquad \downarrow \\
E_{1}\times E_{2} & \stackrel{\varphi_{{O}_{E_1}}\times \varphi_{{O }_{E_2}}}{\longrightarrow}&{\rm Pic}^{0}(E_1) \times {\rm Pic}^{0}(E_2) &
\ni (L\mid_{E_{1}}, L\mid_{E_{2}})
\end{array}
$$
\end{lm}
\proof{
The first statement is well-known (cf. Mumford \cite{M}).
For $x = (x_{1}, x_{2}) \in A$, we have
$$
{\rm Res}\circ \varphi_{X}(x) = {\rm Res}(T_{x}^{*}X - X) \\
= 
(T_{x_{1}}^{*}{O}_{E_1} -{ O}_{E_1}, T_{x_{2}}^{*}{O}_{E_2} -{O}_{E_2})
=(\varphi_{{O}_{E_1}}\times \varphi_{{O}_{E_2}}) (x)
$$
}
We now examine the canonical involution of $B$. Since we have 
$B = {\rm End}(E)\otimes {\bf Q}$, it suffices to define it for the elements
of ${\rm End (E)}$. Then, for $g \in {\rm End} (E)$, the canonical involution is 
given by
$$
    {\bar g} = \varphi^{-1}_{{O}_{E}}\circ {g^{*}}\circ \varphi_{{O}_{E}},
$$
which is the Rosati involution of ${\rm End}(E)\otimes {\bf Q}$ 
(cf. Mumford \cite{M}, Section 21, and Tate \cite{T}, Section 4).
For the elliptic curve $E$, we have
$$
{\bar g}\circ g = \varphi^{-1}_{{O}_{E}}\circ {g^{*}}\circ \varphi_{{O}_{E}}\circ g 
= (\deg g){\rm id}_{E}.
$$
\begin{lm}
Under the Rosati involution the element
$g = 
\left(
\begin{array}{cc}
\alpha & \beta \\
\gamma & \delta 
\end{array}
\right)
\in M_{2}(B)$ maps to
$$
   g' = 
\left(
\begin{array}{cc}
\bar{\alpha} & \bar{\gamma} \\
\bar{\beta} & \bar{\delta} 
\end{array}
\right)
$$
\end{lm}
\proof{We denote by $\hat{g}$ the dual morphism of $g$. 
As the action on divisors, we have $\hat{g} = g^{*}$.
The Rosati involution is given by 
$g' = \varphi^{-1}_{X}\circ \hat{g}\circ \varphi_{X}$.
We calculate the right-hand-side term explicitly. 
We have a  commutative diagram
$$
\begin{array}{ccc}
E_{1} \times E_{2} & 
\stackrel{\varphi_{{O}_{E_1} }\times \varphi_{{O}_{E_2} }}{\longrightarrow} &
{\rm Pic}^{0}(E_{1})\times {\rm Pic}^{0}(E_{2})\\
\quad \downarrow \varphi_{X} &   &   || \\
{\rm Pic}^{0}(E_{1}\times E_{2}) & \stackrel{{\rm Res}}{\longrightarrow} & {\rm Pic}^{0}(E_{1})\times {\rm Pic}^{0}(E_{2}) \\
\uparrow \hat{g} &  & \qquad \qquad \qquad  \uparrow {\rm Res}\circ \hat{g}\circ {\rm Res}^{-1} \\
{\rm Pic}^{0}(E_{1}\times E_{2}) & \stackrel{{\rm Res}}{\longrightarrow} & 
{\rm Pic}^{0} (E_{1})\times {\rm Pic}^{0}(E_{2})  \\
\quad \uparrow \varphi_{X} &  &  || \\
E_{1} \times E_{2} & 
\stackrel{\varphi_{{O}_{E_1} }\times \varphi_{{O}_{E_2} }}{\longrightarrow} 
& {\rm Pic}^{0}(E_{1})\times {\rm Pic}^{0}(E_{2}).
\end{array}
$$
Using this diagram, for the point 
$(x_{1}, x_{2}) \in E_{1}\times E_{2}$
we have
$$
\begin{array}{cl}
g'\left(
\begin{array}{c}
x_{1}\\
x_{2}
\end{array}
\right) & = \varphi^{-1}_{X}\circ \hat{g}\circ \varphi_{X}\left(
\begin{array}{c}
x_{1}\\
x_{2}
\end{array}
\right)  \\
& = (\varphi_{{O}_{E_1} }\times \varphi_{{O}_{E_2} })^{-1}\circ {\rm Res} \circ \hat{g} \circ{\rm Res}^{-1}
\circ (\varphi_{{O}_{E_1} }\times \varphi_{{O}_{E_2} })\left(
\begin{array}{c}
x_{1}\\
x_{2}
\end{array}
\right)  \\
& =(\varphi_{{O}_{E_1} }\times \varphi_{{O}_{E_2} })^{-1}\circ {\rm Res} \circ \hat{g} \circ {\rm Res}^{-1}\left(
\begin{array}{c}
\varphi_{{O}_{E_1}}(x_{1})\\
\varphi_{{O}_{E_2}}(x_{2})
\end{array}
\right)       \\
& = (\varphi_{{O}_{E_1} }\times \varphi_{{ O}_{E_2} })^{-1}\circ {\rm Res} \circ \hat{g}(p_{1}^{*}\varphi_{{\cal O}_{1}}(x_{1}) +  p_{2}^{*}\varphi_{{O}_{E_2}}(x_{2})) \\
& = (\varphi_{{O}_{E_1} }\times \varphi_{{O}_{E_2} })^{-1}\circ {\rm Res}((p_{1}\circ g)^{*}\varphi_{{O}_{E_1}}(x_{1}) +  (p_{2}\circ g)^{*}\varphi_{{O}_{E_2}}(x_{2})) 
\end{array}
$$
We denote by $m_{i}$ the addition of $E_{i}$ $(i = 1, 2)$. Then, we have
$$
 p_{1}\circ g = m_{1}\circ (\alpha \times \beta),~ p_{2}\circ g = m_{2}\circ (\gamma \times \delta).
$$
We denote by $q_{i}$ $(i = 1, 2)$ the i-th projection $E_{1}\times E_{1}\rightarrow E_{1}$. Then by Mumford \cite{M}, for $L\in {\rm Pic}^{0}(E_{1})$ we have
$$
   m_{1}^{*}L \sim q_{1}^{*}L + q_{2}^{*}L~\mbox{(linearly equivalent)}.
$$
Therefore we have
$$
\begin{array}{cl}
(p_{1}\circ g)^{*}\varphi_{{O}_{E_1}}(x_{1}) &  = (m_{1}\circ (\alpha \times \beta))^{*}\varphi_{{O}_{E_1}}(x_{1}) \\
  & = (\alpha \times \beta)^{*}m_{1}^{*}\varphi_{{O}_{E_1}}(x_{1})\\
& = (\alpha \times \beta)^{*}(q_{1}^{*}\varphi_{{O}_{E_1}}(x_{1}) +
 q_{2}^{*}\varphi_{{O}_{E_1}}(x_{1})) \\
&=  \{q_{1}\circ (\alpha \times \beta)\}^{*}\varphi_{{O}_{E_1}}(x_{1}) +
 \{q_{2}\circ (\alpha \times \beta)\}^{*}\varphi_{{O}_{E_1}}(x_{1}).
\end{array}
$$
Since we have commutative diagrams
$$
\begin{array}{ccc}
E_{1}\times E_{2} & \stackrel{\alpha\times \beta}{\longrightarrow} & 
E_{1}\times E_{1} \\
\downarrow p_{1} &    & \downarrow q_{1} \\
E_{1} &  \stackrel{\alpha}{\longrightarrow} & E_{1},
\end{array}
\quad
\begin{array}{ccc}
E_{1}\times E_{2} & \stackrel{\alpha\times \beta}{\longrightarrow} & 
E_{1}\times E_{1} \\
\downarrow p_{2} &    & \downarrow q_{2} \\
E_{1} & \stackrel{\beta}{\longrightarrow} & E_{1},
\end{array}
$$
we have
$$
(p_{1}\circ g)^{*}\varphi_{{O}_{E_1}(x_{1})}  = p_{1}^{*}\alpha^{*}\varphi_{{O}_{E_1}(x_{1})} + p_{2}^{*}\beta^{*}\varphi_{{ O}_{E_1}(x_{1})}.
$$
In a similar way, we have
$$
(p_{2}\circ g)^{*}\varphi_{{O}_{E_2}(x_{2})}  = p_{1}^{*}\gamma^{*}\varphi_{{ O}_{E_2}(x_{2})} + p_{2}^{*}\delta^{*}\varphi_{{O}_{E_2}(x_{2})}.
$$
Therefore, we have
$$
\begin{array}{cl}
g'\left(
\begin{array}{c}
x_{1}\\
x_{2}
\end{array}
\right) &
=(\varphi_{{O}_{E_1} }\times \varphi_{{O}_{E_2} })^{-1}
  \circ {\rm Res}
(p_{1}^{*}\alpha^{*}\varphi_{{O}_{E_1}}(x_{1}) + p_{2}^{*}\beta^{*}\varphi_{{ O}_{E_1}}(x_{1}) \\
 & \qquad  \qquad+ p_{1}^{*}\gamma^{*}\varphi_{{O}_{E_2}}(x_{2}) + p_{2}^{*}\delta^{*}\varphi_{{O}_{E_2}}(x_{2}))\\
 & = (\varphi_{{O}_{E_1} }\times \varphi_{{O}_{E_2}})^{-1}
\left(
\begin{array}{c}
\alpha^{*}\varphi_{{O}_{E_1}}(x_{1}) + 
\gamma^{*}\varphi_{{O}_{E_2}}(x_{2})\\
 \beta^{*}\varphi_{{O}_{E_1}}(x_{1})+ 
\delta^{*}\varphi_{{O}_{E_2}}(x_{2})
\end{array}
\right) \\
&= \left(
\begin{array}{c}
\varphi_{{O}_{E_1} }^{-1}\alpha^{*}\varphi_{{O}_{E_1}}(x_{1}) + 
\varphi_{{O}_{E_1} }^{-1}\gamma^{*}\varphi_{{O}_{E_2}}(x_{2})\\
\varphi_{{O}_{E_2} }^{-1}\beta^{*}\varphi_{{O}_{E_1}}(x_{1}) + 
\varphi_{{O}_{E_2} }^{-1}\delta^{*}\varphi_{{O}_{E_2}}(x_{2})
\end{array}
\right) \\
& = \left(
\begin{array}{cc}
\varphi_{{O}_{E_1} }^{-1}\alpha^{*}\varphi_{{O}_{E_1}} & 
\varphi_{{O}_{E_1} }^{-1}\gamma^{*}\varphi_{{O}_{E_2}}\\
\varphi_{{O}_{E_2} }^{-1}\beta^{*}\varphi_{{O}_{E_1}} &
\varphi_{{O}_{E_2} }^{-1}\delta^{*}\varphi_{{O}_{E_2}}
\end{array}
\right) 
\left(
\begin{array}{c}
x_{1}\\
x_{2}
\end{array}
\right) 
\end{array}
$$
Since $E_{1} = E_{2} = E$ and $\varphi_{{\cal O}_{1}}= \varphi_{{\cal O}_{2}}$,
we conclude
$$
   g'= 
\left(
\begin{array}{cc}
{\alpha} & {\beta} \\
{\gamma} & {\delta} 
\end{array}
\right)'
=
\left(
\begin{array}{cc}
\bar{\alpha} & \bar{\gamma} \\
\bar{\beta} & \bar{\delta} 
\end{array}
\right)
$$}
\begin{lm}\label{lm;E}
For a divisor $L \in {\rm Pic}(E_{1}\times E_{2})$ with  $j(L) = g = 
\left(
\begin{array}{cc}
{\alpha} & {\beta} \\
{\gamma} & {\delta} 
\end{array}
\right)
$,
we have 
$$
\alpha = L\cdot E_{1}, \quad  \delta = L\cdot E_{2}.
$$
\end{lm}
\proof{
Since $\alpha$ is an integer, we have
$$
(1) \qquad g\left(
\begin{array}{c}
x\\
{O}_{E_2}
\end{array}
\right)
= 
\left(
\begin{array}{cc}
{\alpha} & {\beta} \\
{\gamma} & {\delta} 
\end{array}
\right)
\left(
\begin{array}{c}
x\\
{O}_{E_2}
\end{array}
\right)
=
\left(
\begin{array}{c}
\alpha x\\
\gamma(x)
\end{array}
\right).
$$
Now, we examine $\alpha x$. 
$$
\begin{array}{cl}
   g\left(
\begin{array}{c}
x\\
{O}_{E_2}
\end{array}
\right) 
& = \varphi^{-1}_{X}\circ \varphi_{L}\left(
\begin{array}{c}
x\\
{O}_{E_2}
\end{array}
\right) \\
& = \varphi^{-1}_{X}\{T^{*}_{(x, {\cal O}_{2})}L - L\} \\
& = (\varphi_{{\cal O}_{1} }\times \varphi_{{\cal O}_{2} })^{-1}\circ 
{\rm Res}\{T^{*}_{(x, {\cal O}_{2})}L - L\}
\end{array}
$$
We restrict the divisor $L$ to $E_{1}$ and denote it by $e$. Then, the divisor is
expressed as
$$
     e \sim \sum_{i=1}^{\lambda}n_{i}P_{i}
$$
with integers $n_{i}$ and points $P_{i}$ on $E_{1}$ $(i = 1, 2, \cdots, \lambda)$.
We have
$$
L\cdot E_{1} = \deg e = \sum_{i = 1}^{\lambda}n_{i}.
$$
We set $n = \sum_{i = 1}^{\lambda}n_{i}$. 
Then, we obtain the following form:
$$
   g\left(
\begin{array}{c}
x\\
{O}_{E_2}
\end{array}
\right)
= (\varphi_{{O}_{E_1} }\times \varphi_{{O}_{E_2} })^{-1}
\left(
\begin{array}{c}
T^{*}_{x}e - e\\
  {}*
\end{array}
\right).
$$
We denote by $\oplus$ the addition of $E_{1}$, and by $\ominus$ 
the subtraction of $E_{1}$. Then, we have
$$
   T_{x}^{*}e \sim \sum_{i = 1}^{\lambda}n_{i}(P_{i} \ominus x).
$$
By Abel's theorem, we see that
$$
\begin{array}{cl}
T_{x}^{*} e - e  &
\sim n_{1}(P_{1}\ominus x)\oplus \cdots \oplus n_{\lambda}(P_{\lambda}\ominus x)
\ominus (n_{1}P_{1} \oplus \cdots \oplus n_{\lambda}P_{\lambda}) - {O}_{E_1} \\
 & \sim (-n)x - {O}_{E_1} \\
 & = T^{*}_{nx}{O}_{E_1} - {O}_{E_1} \\
 & = \varphi_{{O}_{E_1}}(nx).
\end{array}
$$
Therefore, we have  $\varphi_{{O}_{E_1}}^{-1}((-nx) -{O}_{E_1}) = nx$,
and 
$$
(2) \qquad    g\left(
\begin{array}{c}
x\\
{O}_{E_2}
\end{array}
\right) =  
\left(
\begin{array}{c}
nx \\
{}*
\end{array}
\right) . 
$$
Hence, comparing (1) and (2), we have $\alpha = n = L\cdot E_{1}$.
In a similar way, we have $\delta = L\cdot E_{2}$.
}
\begin{lm}
We have $j (E_{1}) = \left(
\begin{array}{cc}
0 & 0 \\
0 & 1
\end{array}
\right)$ and 
$j (E_{2}) = \left(
\begin{array}{cc}
1 & 0 \\
0 & 0
\end{array}
\right)$
\end{lm}
\proof{
For a point $(x_{1}, x_{2}) \in E_{1}\times E_{2}$, we have
$$
\begin{array}{cl}
\varphi_{X}^{-1}\circ \varphi_{E_{1}}
\left(
\begin{array}{c}
x_{1}\\
x_{2}
\end{array}
\right)
& = \varphi_{X}^{-1}\{T^{*}_{(x_{1}, x_{2})}E_{1} - E_{1}\}\\
& = (\varphi_{{O}_{E_1} }\times \varphi_{{O}_{E_2} })^{-1}\circ 
{\rm Res}\{T^{*}_{(x_{1}, x_{2})}E_{1} - E_{1}\}\\
& = (\varphi_{{O}_{E_1} }\times \varphi_{{O}_{E_2} })^{-1}
\left(
\begin{array}{c}
{O}_{E_1} - {O}_{E_1}\\
(-x_{2}) -{O}_{E_2}
\end{array}
\right)
 = \left(
\begin{array}{c}
{O}_{E_1}\\
x_{2}
\end{array}
\right).
\end{array}
$$
Therefore, we have
$$
j(E_{1}) =
\varphi_{X}^{-1}\circ \varphi_{E_{1}} = \left(
\begin{array}{cc}
0 & 0 \\
0 & 1
\end{array}
\right).
$$
In a similar way, we obtain the second assertion.
}
\begin{lm}
For $L \in {\rm NS}(E_{1}\times E_{2})$, we set $j(L) = g = \left(
\begin{array}{cc}
\alpha & \bar{\gamma} \\
\gamma & \delta
\end{array}
\right)$.
 Then,
$$
     L^{2} = 2 \det g.
$$
\end{lm}
\proof{
Since $\alpha, \delta \in {\bf Z}$, we have
$$
\begin{array}{cl}
\varphi^{-1}_{X}\circ \varphi_{(L -\alpha E_{2} -\delta E_{1})}
& = \varphi^{-1}_{X}\circ \varphi_{L} - \alpha\varphi^{-1}_{X}\circ \varphi_{E_{2}}
        -\delta \varphi^{-1}_{X}\circ\varphi_{E_{1}} \\
  &= \left(
\begin{array}{cc}
\alpha & \bar{\gamma} \\
\gamma & \delta
\end{array}
\right)
- \left(
\begin{array}{cc}
\alpha & 0 \\
0 & 0
\end{array}
\right)
- \left(
\begin{array}{cc}
0 & 0 \\
0 & \delta
\end{array}
\right) \\
& =  \left(
\begin{array}{cc}
0 & \bar{\gamma} \\
\gamma & 0
\end{array}
\right).
\end{array}
$$
Since the right hand-side is contained in $H$, there exists a divisor $Z$
such that 
$$
\varphi^{-1}_{X}\circ \varphi_{Z} = \left(
\begin{array}{cc}
0 & \bar{\gamma} \\
\gamma & 0
\end{array}
\right).
$$
If $Z$ is zero, then we have $\gamma = 0$. Therefore, we have 
$ Z = \alpha E_{2} +\delta E_{1}$ and $Z^2 = 2\alpha\delta = \det g$.

Now, we assume $Z\neq 0$.
Since $\varphi_{X}$ is an isomorphism, by the Riemann-Roch theorem on the abelian surface $A$,
we have
$$
\deg (\varphi^{-1}_{X}\circ \varphi_{Z}) = \deg \varphi_{Z} =(Z^{2}/2)^{2}.
$$
On the other hand,
$$
\deg (\varphi^{-1}_{X}\circ \varphi_{Z}) =\deg\gamma \cdot \deg\bar{\gamma} =
(\deg \gamma)^{2} = (\gamma\bar{\gamma})^2 
$$
By Lemma \ref{lm;E}, we have
$$
Z\cdot E_{1} = Z\cdot E_{2} = 0.
$$
Therefore, we have $Z\cdot (E_{1} + E_{2}) = 0$. Since $(E_{1} + E_{2})^{2} = 2 > 0$,
by the Hodge index theorem we see $Z^{2} < 0$. Therefore, we have,
$Z^{2}/2 = -\gamma\bar{\gamma}$.

On the other hand, since $\varphi_{X}$ is an isomorphism and $\varphi^{-1}_{X}\circ \varphi_{(L -\alpha E_{2} -\delta E_{1} -Z)} = 0$, we have $\varphi_{(L -\alpha E_{2} -\delta E_{1} -Z)} = 0$. Therefore, we have 
$$
     0 \equiv L -\alpha E_{2} -\delta E_{1} -Z,
$$
where by $\equiv$ we mean algebraic equivalence.
Hence, we have
$$
L^{2} = 2\alpha \delta + Z^{2} = 2(\alpha \delta - \gamma\bar{\gamma}) = 2\det g.
$$
}
For an automorphism $g$ of $A$, we can regard $g$ as an element of
${\rm M}_{2}({\cal O})$, and then we can consider ${}^{t}\bar{g}$.
\begin{lm}
Let $L_{1}$ and $L_{2}$ be two divisors with $j(L_{1}) = g_{1}$ and
$j(L_{2}) = g_{2}$. Let $g$ be an automorphism of $A$.
Then, $g^{*}L_{1} \equiv L_{2}$ if and only if ${}^{t}\bar{g} g_{1}g = g_{2}$.
\end{lm}
\proof{
We have 
$$
\begin{array}{ccl} 
g^{*}L_{1} \equiv L_{2} &  \Longleftrightarrow & \varphi_{g^{*}L_{1}} = \varphi_{L_{2}}\\
& \Longleftrightarrow & \hat{g}\circ \varphi_{L_{1}}\circ g = \varphi_{L_{2}}\\
& \Longleftrightarrow & \varphi^{-1}_{X} \circ \hat{g}\circ \varphi_{X}\circ (\varphi^{-1}_{X}\circ \varphi_{L_{1}})\circ g =\varphi_{X}^{-1} \circ \varphi_{L_{2}} \\
& \Longleftrightarrow &g'\circ g_{1}\circ g = g_{2}.
\end{array}
$$
}
Let  $m : E\times E \rightarrow E$ be the addition of $E$, and we set
$$
\Delta = {\rm Ker}~ m.
$$
We have $\Delta = \{(P, - P)  ~\mid~P\in E \}$. Note that this $\Delta$ is 
different from the usual diagonal. For two endomorphisms
$a_{1}, a_{2}\in {\rm End}(E)$, we set
$$
\Delta_{a_{1}, a_{2}} = (a_{1}\times a_{2})^{*}\Delta.
$$
Using this notation, we have $\Delta = \Delta_{1,1}$.
We have the following theorem (cf. \cite{K2}).
\begin{th}\label{div}
$$
j(\Delta_{a_{1}, a_{2}}) =
\left(
\begin{array}{cc}
\bar{a}_{1}a_{1}  & \bar{a}_{1}a_{2}\\
\bar{a}_{2}a_{1} & \bar{a}_{2}a_{2}
\end{array}
\right).
$$
In particular, we have
$$
j(\Delta) =
\left(
\begin{array}{cc}
1  & 1\\
1 & 1
\end{array}
\right).
$$
\end{th}
\proof{
Let $\alpha, \beta, \gamma$ be elements of ${\cal O}$ such that 
$$
\varphi_{X}^{-1}\circ \varphi_{\Delta} = \left(
\begin{array}{cc}
\alpha & \bar{\gamma} \\
\gamma & \delta
\end{array}
\right).
$$
Then, since $E_{1}\cdot \Delta = E_{2}\cdot \Delta = 1$, we have
$\alpha = \delta = 1$ by Lemma \ref{lm;E}. Since we have
$$
\varphi_{X}^{-1}\circ \varphi_{\Delta}
\left(
\begin{array}{c}
x\\
-x
\end{array}
\right) 
= \varphi^{-1}_{X}\{T^{*}_{(x, -x)}\Delta - \Delta\} = \varphi_{X}^{-1}(0) =
\left(
\begin{array}{c}
{\cal O}_{1}\\
{\cal O}_{2}
\end{array}
\right), 
$$
we have $\gamma (x) = x$ for any $x \in E$.
Therefore, we have $\gamma = 1$.

By definition, we have
$$
\Delta_{a_{1}, a_{2}} =(a_{1}\times a_{2})^{*}\Delta.
$$
Therefore, we have
$$
j(\Delta_{a_{1}, a_{2}}) = 
{}^{t}\overline{\left(
\begin{array}{cc}
a_{1}  & 0\\
0 & a_{2}
\end{array}
\right)}
\left(
\begin{array}{cc}
1  & 1\\
1 & 1
\end{array}
\right)
\left(
\begin{array}{cc}
a_{1}  & 0\\
0 & a_{2}
\end{array}
\right)
=
\left(
\begin{array}{cc}
\bar{a}_{1}a_{1}  & \bar{a}_{1}a_{2}\\
\bar{a}_{2}a_{1} & \bar{a}_{2}a_{2}
\end{array}
\right)
$$
}
\section{Non-superspecial cases}
In this section, we examine the injectivity of the Chern class map 
of abelian surfaces. Let $\alpha_{p}$ be the local-local group scheme
of rank $p$ (cf. Oort \cite{Oo1} for the definition and properties). Then, we have
${\rm End}(\alpha_{p}) \simeq k$, and for an abelian variety $X$, 
${\rm Hom}(\alpha_{p}, X)$ is a right vector space over 
${\rm End}(\alpha_{p})\simeq k$ by composition of morphisms.
The a-number of $X$ is defined by
$$
        a = \dim_{k}{\rm Hom}(\alpha_{p}, X).
$$
We denote by $[p]_{X}$ multiplication of $p$:
$$
\begin{array}{rccc}
    [p]_{X} : & X  & \longrightarrow & X \\
       &   x & \mapsto & px.
\end{array}
$$
Then, the reduced part of $\Ker ~[p]_{X}$ is of the form:
$$
(\Ker ~[p]_{x})_{red} \simeq ({\bf Z}/p{\bf Z})^{\oplus r}
$$
with an integer $r$ $(0 \leq r \leq \dim X$). We call $r$ the p-rank of $X$
(cf. Mumford \cite{M}).
The following theorem follows essentially from the results in Ogus \cite{O}, 
but we give here a down-to-earth proof. For the definition and properties
of the Cartier operator, see Cartier \cite{C}.
\begin{th}\label{injectivity} 
Let $X$ be an abelian surface defined over $k$. Then,
the Chern class map
$$
     c_{1}: {\rm NS}(X)/p{\rm NS}(X) \rightarrow {\rm H}^{1}(X, \Omega_{X}^1)
$$
is injective if and only if $X$ is not superspecial.
\end{th}
\proof{The only-if-part will be proved in Theorem \ref{mainth}. 
We prove here the if-part.
We denote by
$r(X)$ the p-rank of $X$, and by $a(X)$ the a-number of $X$. 
By Oort \cite{Oo2}, $X$ is superspecial if and only if $a(X) = 2$.
Therefore, we assume $a(X) \neq 2$.
Take an affine open covering $\{U_{i}\}$ of $X$, and suppose that there is
a divisor $D = \{f_{ij}\}$ which is not zero in ${\rm NS}(X)/p{\rm NS}(X)$, such that
$c_{1}(D) = \{df_{ij}/f_{ij}\} \sim 0$ in ${\rm H}^{1}(X, \Omega_{X}^1)$. 
Then, there exists $\omega_{i} \in  {\rm H}^{0}(U_{i}, \Omega_{X}^1)$
such that
$$
            df_{ij}/f_{ij} = \omega_{j} - \omega_{i}.
$$

(i) The first case : $d\omega_{i} = 0$.

 Applying the Cartier operator $C$, we obtain
$$
            df_{ij}/f_{ij} = C(\omega_{j}) - C(\omega_{i}).
$$
Therefore, we have
$$
       C(\omega_{j}) - \omega_{j} = C(\omega_{i}) -\omega_{i} \quad
\mbox{on}~U_{i}\cap U_{j}.
$$
and we have a regular 1-form $\omega'$ on $X$ which is defined by
$$
C(\omega_{i}) -\omega_{i} \quad \mbox{on}~U_{i}.
$$
Since $C - {\rm id}: H^{0}(X, \Omega^{1}_{X}) \rightarrow H^{0}(X, \Omega^{1}_{X})$
is surjective, there exists a regular 1-form $\omega \in H^{0}(X, \Omega^{1}_{X})$
such that $(C - {\rm id})(\omega) = \omega'$. Therefore, we have
$$
      C(\omega_{i}- \omega) = \omega_{i} - \omega.
$$
By the property of the Cartier operator, there exists an regular function 
$f_{i}$ on $U_{i}$ such that
$$
     \omega_{i}- \omega = df_{i}/f_{i},
$$
and we have 
$$
      df_{ij}/f_{ij} =   df_{j}/f_{j} - df_{i}/f_{i}.
$$
This means $d(f_{ij}f_{i}/f_{j}) = 0$. Therefore, there exists a regular function
$g_{ij}$ on $U_{i}\cap U_{j}$ such that
$$
       f_{ij}f_{i}/f_{j}  = g_{ij}^p   \quad \mbox{on}~U_{i}\cap U_{j}.
$$
Since $D = \{  f_{ij}\}$ is a cocycle, we see that $\{g_{ij} \}$ is also a cocyle
and that this gives an element of ${\rm NS}(X)$.
Therefore, we conclude $D \in p{\rm NS}(X)$, which
contradicts $D \neq 0$ in ${\rm NS}(X)/p{\rm NS}(X)$.

(ii) The second case : $d\omega_{i} \neq 0$.

In this case we have $d\omega_{i} = d\omega_{j}$ on $U_{i}\cap U_{j}$
and we get a non-zero regular 2-form on $X$.
Since this regular 2-form is d-exact and is a basis of 
${\rm H}^{0}(X, \Omega_{X}^{2})$, the Cartier operator acts
on ${\rm H}^{0}(X, \Omega_{X}^{2})$ as the zero map.
Therefore, $X$ is not ordinary, that is, $r(X) \neq 2$.
Therefore, we have either $r(X) = 1$ and $a(X) = 1$,
or $r(X) = 0$ and $a(X) = 1$. 

Now, we consider the absolute Frobenius
$F: {\rm H}^{1}(X, {\cal O}_{X}) \rightarrow
{\rm H}^{1}(X, {\cal O}_{X})$. Since $a(X) = 1$ in both cases, there exists 
a non-zero element $\beta= \{g_{ij}\}$ in ${\rm H}^{1}(X, {\cal O}_{X})$ 
such that $F(\beta) = 0$. 
This means that there exists a regular function $g_{i}$
on $U_{i}$ such that $g_{ij}^p = g_{j} - g_{i}$. Since $dg_{i} = dg_{j}$ 
on $U_{i} \cap U_{j}$, we have a non-zero regular 1-form $\eta$ on $X$
given by $dg_{i}$ on $U_{i}$.  Since $\dim {\rm H}^{0}(X, \Omega_{X}^{1})  = 2$,
in both cases there exists a nonzero regular 1-form $\eta'$
such that  $\{\eta, \eta'\}$ gives a basis of ${\rm H}^{0}(X, \Omega_{X}^{1})$
with $C(\eta') \neq 0$. In fact, we can take $\eta'$ with $C(\eta') = \eta$
if $r(X) =0$  and $a(X) = 1$, and we can take $\eta'$ with $C(\eta') = \eta'$
if $r(X) =1$ and $a(X) = 1$. Since we have 
${\rm H}^{0}(X, \Omega_{X}^{2}) = \wedge^{2}{\rm H}^{0}(X, \Omega_{X}^{1}) $,
$\eta\wedge \eta'$ gives a basis of ${\rm H}^{0}(X, \Omega_{X}^{2})$.
Therefore, there exists a non-zero element $a \in k$ such that 
$$
          d\omega_{i} = a\eta \wedge \eta' = a(d(g_{i}\eta')).
$$
We set $\theta_{i} = \omega_{i} - ag_{i}\eta'$. Then, $\theta_{i}$ is
d-closed and we have
$$
\begin{array}{cl}
 df_{ij}/f_{ij} & =  ag_{j}\eta' -ag_{i}\eta' + \theta_{j} - \theta_{i}\\
                      & = ag_{ij}^{p}\eta' + \theta_{j} - \theta_{i}
\end{array}
$$
Applying the Cartier operator, we have
$$
df_{ij}/f_{ij} = a^{1/p}g_{ij}C(\eta') + C(\theta_{j}) - C(\theta_{i}).
$$
This means that 
$$
c_{1}(D) \sim a^{1/p}\beta\otimes C(\eta') \in 
{\rm H}^{1}(X, {\cal O}_{X})\otimes {\rm H}^{0}(X, \Omega_{X}^{1})\cong
{\rm H}^{1}(X, \Omega_{X}^1)
$$
Since $\beta \neq 0$ in ${\rm H}^{1}(X, {\cal O}_{X})$
and $C(\eta') \neq 0$ in ${\rm H}^{0}(X, \Omega^1_{X})$, 
we see  $\beta\otimes C(\eta')\neq 0$
in ${\rm H}^{1}(X, \Omega_{X}^1)$. A contradiction.

Hence, if $a(X) \neq 2$, we conclude that $c_{1}$ is injective.
}

\section{Superspecial cases}
Let $k$ be an algebraically closed field of characteristic $p \geq 3$.
For an elliptic curve $E$ over $k$,
we examine the action of endomorphisms of $E$ on 
${\rm H}^{0}(E, \Omega^{1}_{E})$ and ${\rm H}^{1}(E, {\cal O}_{E})$.
\begin{lm}\label{3.1}
Let $E$ be an elliptic curve and $\alpha \in {\rm End}(E)$.
Assume $\alpha$ acts on ${\rm H}^{1}(E, {\cal O}_{E})$ as multiplication by 
$\beta \in k$ 
$(\beta \neq  0)$.
Then, $\alpha$ acts on ${\rm H}^{0}(E, \Omega^{1}_{E})$ as multiplication by ${\deg \alpha}/{\beta}$.
\end{lm}
\proof{Using the endomorphism $\alpha : E \rightarrow E$, we obtain
a commutative diagram
$$
\begin{array}{ccc}
{\rm NS}(E)/p{\rm NS}(E) & \stackrel{\alpha^*}{\longrightarrow}& 
{\rm NS}(E)/p{\rm NS}(E) \\
  \downarrow c_{1}  &        & \downarrow c_{1}\\
{\rm H}^{1}(E, \Omega^{1}_{E}) &   \stackrel{\alpha^*}{\longrightarrow}& 
{\rm H}^{1}(E, \Omega^{1}_{E})\\
  \downarrow &        & \downarrow \\
{\rm H}^{1}(E, {\cal O}_{E}) \otimes {\rm H}^{0}(E, \Omega^{1}_{E})& \stackrel{\alpha^{*}  \otimes \alpha^{*}}{\longrightarrow}& {\rm H}^{1}(E, {\cal O}_{E}) \otimes {\rm H}^{0}(E, \Omega^{1}_{E}).
\end{array}
$$
Take a point $Q \in E$, and bases $\omega \in {\rm H}^{0}(E, \Omega^{1}_{E})$,
$\eta \in {\rm H}^{1}(E, {\cal O}_{E})$. Then, we have 
$ \alpha^{*}(Q)=(\deg \alpha)Q $,
and $ (\alpha^{*} \otimes \alpha^{*})(\omega \otimes \eta)=(\beta \omega) \otimes \alpha^{*}\eta$. The result follows from the diagram.
}
For an integer $n$, we have an endomorphism
$[n]_E : E \longrightarrow E$ given by $P \mapsto nP$ $(P \in E)$.
\begin{lm}The induced homomorphism 
$$
[n]^*_E: {\rm H}^{0}(E, \Omega^{1}_{E}) \longrightarrow {\rm H}^{0}(E, \Omega^{1}_{E})
$$ 
is multiplication by $n$, 
i.e., $[n]^{*}_{E}\omega=n\omega$ for $\omega \in H^{0}(E, \Omega^{1}_{E})$.
\end{lm}
\proof{This follows from the fact that $[n]_{*}$ is given as multiplication by 
$n$ on the tangent space at the origin (Mumford \cite{M}).
}

Assume $p \ne 2$. Following the theory of Ibukiyama (cf. \cite{I}) to construct
a quaternion division algebra over ${\bf Q}$ with discriminant $p$,
we take a prime number $q$ such that 
$-q \equiv 5~({\rm mod}~ 8)$ and $(\frac{-q}{p})=-1$, and take an integer $a$ 
such that $a^{2} \equiv-p~({\rm mod}~q)$. 
Here, $(\frac{-q}{p})$ is the Legendre symbol.
Then, the quaternion division algebra
$B$ over ${\bf Q}$ with discriminant $p$ and a maximal order ${\cal O}$ of $B$
are given by
$$
\begin{array}{c}
B={\bf Q} \oplus {\bf Q}F \oplus {\bf Q} \alpha \oplus {\bf Q}F \alpha  \\
\mbox{with}~ F^{2}=-p, \alpha^{2}=-q, F \alpha =-\alpha F \\
{\cal O} = {\bf Z}+{\bf Z}(\frac{1+\alpha}{2})+{\bf Z}(\frac{F(1+\alpha)}{2})
+{\bf Z}(\frac{(a+F)\alpha}{q}).
\end{array}
$$
Then, we know that there exists a supersingular elliptic curve 
$E$ over $k$ with ${\rm End}(E) = {\cal O}$ and ${\rm End}^0(E) = B$
(cf. Deuring \cite{D}). 

We need the following  well-known lemma. 
\begin{lm}\label{lm,chern}
For a non-singular complete algebraic curve X,  the Chern class map
$$
c_{1} : {\rm Pic}(X)/p{\rm Pic}(X) \hookrightarrow   {\rm H}^{1}(X, \Omega^{1}_{X})
$$
is injective.
\end{lm}
\proof{Let $L$ be a class of  ${\rm Pic}(X)/p{\rm Pic}(X)$.
Then, we can lift this class to ${\rm Pic}(X)$. We take an open affine
covering $\{U_{i}\}$ that trivializes the corresponding invertible sheaf, and
let the invertible sheaf be
given by $\{f_{ij}\}$ with a regular function $f_{ij}$ on $U_{i}\cap U_{j}$.
Then, we have $c_{1}(L) = \{df_{ij}/f_{ij}\}$.

Suppose $\{df_{ij}/f_{ij}\} \sim 0$. 
Then, there exists ${\omega_{i}\in \Omega^{1}_{X}(U_{i})}$
such that 
$$ 
\frac{df_{ij}}{f_{ij}}=\omega_{j}-\omega_{i}. 
$$ 
Since $X$ is one-dimensional,
$\omega_{i}$'s are d-closed. By the Cartier operator $C$, we have
$$
\frac{df_{ij}}{f_{ij}} = C(\omega_{j}) - C(\omega_{i}).
$$
Therefore, we have 
$$
C(\omega_{i})-\omega_{i}=C(\omega_{j})-\omega_{j}.
$$ 
Hence, $C(\omega_{i})-\omega_{i}$ on $U_{i}$ gives
a global regular 1-form $\omega \in {\rm H}^{0}(X,  \Omega^{1}_{X})$.
Since  $C - {\rm id}_{X}$ is surjective on ${\rm H}^{0}(X, \Omega^{1}_{X})$, 
there exists
$\{\tilde{\omega}\} $ such that $(C - {\rm  id}_{X})(\tilde{\omega})=\omega$.
Replace $\omega_{i}$
by $\omega_{i}-\tilde{\omega}$, we may assume $C(\omega_{i}) = \omega_{i}$.
Hence, there exists
$f_{i}$ such that $\omega_{i}=\frac{df_{i}}{f_{i}}$. 
The result follows from this fact (cf. the proof of Theorem \ref{injectivity}).
}

We now compute the Chern class map explicitly for $A$, 
where $A=E_{1}\times E_{2}$ with $E_{1}=E_{2}=E$, 
the supersingular elliptic curve.
The cup product induces a natural isomorphism
$$ 
 {\rm H}^{1}(A, \Omega^{1}_{A}) \cong {\rm H}^{1}(A, {\cal O}_{A}) \otimes 
{\rm H}^{0}(A, \Omega^{1}_{A})
$$
with
$$
\begin{array}{l}
{\rm H}^{1}(A, {\cal O}_{A}) \cong {\rm H}^{1}(E_{1}, {\cal O}_{E_1}) \oplus {\rm H}^{1}(E_{2}, {\cal O}_{E_2}), \\
{\rm H}^{0}(A, \Omega^{1}_{A})\cong {\rm H}^{0}(E_{1}, \Omega^{1}_{E_1}) \oplus {\rm H}^{0}(E_{2}, \Omega^{1}_{E_2}).
\end{array}
$$
Therefore, we have a decomposition
$$
\begin{array}{cl}
(*)\quad {\rm H}^{1}(A, \Omega^{1}_{A})  \cong &({\rm H}^{1}(E_1, {\cal O}_{E_1}) \otimes {\rm H}^{0}(E_1, {\cal O}_{E_1})) \oplus ({\rm H}^{1}(E_1, {\cal O}_{E_1}) \otimes {\rm H}^{0}(E_2, \Omega^{1}_{E_2}))  \\
  &  \oplus  ({\rm H}^{1}(E_2, {\cal O}_{E_2})  \otimes  {\rm H}^{0}(E_1, \Omega^{1}_{E_1}))  \oplus ({\rm H}^{1}(E_2, {\cal O}_{E_2}) \otimes {\rm H}^{0}(E_2, \Omega^{1}_{E_2})).
\end{array}
$$
We have projections 
$$
pr_{i}: A \longrightarrow E_{i}\quad (i = 1, 2).
$$
Then, we have injective homomorphisms 
$$
pr_{i}^{*}:{\rm H}^{1}(E_{i}, \Omega^{1}_{E_i}) 
\hookrightarrow {\rm H}^{1}(A, \Omega^{1}_{A}).
$$
Note that
$$
\begin{array}{l}
     {\rm H}^{1}(E_{1}, \Omega_{E_{1}}^{1}) \cong {\rm H}^{1}(E_{1}, {\cal O}_{E_1})
\otimes {\rm H}^{0}(E_{1}, \Omega^{1}_{E_1})\\
 {\rm H}^{1}(E_{2}, \Omega_{E_{2}}^{1}) \cong {\rm H}^{1}(E_{2}, {\cal O}_{E_2})
\otimes {\rm H}^{0}(E_{2}, \Omega^{1}_{E_2}),
\end{array}
$$
and we have the following commutative diagram
$$
\begin{array}{rccc}
(**)\quad \quad \quad &{\rm NS}(A)/p{\rm NS}(A)& \stackrel{c_{1}}{\hookrightarrow} & {\rm H}^{1}(A, \Omega^{1}_{A})\\
 & \uparrow &                   &pr_{i}^{*}\uparrow            \\
 &{\rm Pic}({E}_{i})/p{\rm Pic}({E}_{i})&\stackrel{c_{1}}{\hookrightarrow} &{\rm H}^{1}({E}_{i}, \Omega^{1}_{{E}_{i}})
\end{array}
$$
The image of the homomorphism $pr_{i}^{*}$ is a one-dimensional subspace 
${\rm H}^{1}(E_{i}, {\cal O}_{E_i}) \otimes {\rm H}^{0}(E_{i}, \Omega^{1}_{E_i})$ 
$(i = 1, 2)$ in ${\rm H}^{1}(A, \Omega^{1}_{A})$.

Now, we consider the Chern class map
$$
 {\rm NS}(A)/p{\rm NS}(A)\cong {\rm Pic}(A)/p{\rm Pic}(A) \stackrel{c_{1}}{\longrightarrow} {\rm H}^{1}(A, \Omega^{1}_{A}).
$$
For the divisors $E_{2}$ (resp. $E_{1}$) on $A$, we set $\Omega_1= c_{1}(E_{2})$ (resp. $\Omega_4 = c_{1}(E_{1})$). Then, by the diagram $(**)$
$\Omega_1$ (resp. $\Omega_4$)
is a basis of ${\rm H}^{1}(E_1, {\cal O}_{E_1}) \otimes {\rm H}^{0}(E_1, {\cal O}_{E_1})$ (resp. ${\rm H}^{1}(E_{2}, {\cal O}_{E_2})
\otimes {\rm H}^{0}(E_{2}, \Omega^{1}_{E_2})$).
 
We set
$$ 
\Delta_a= \Delta_{{\rm id}, a}.
$$
Here,  ${\rm id}$  is the identity endomomorphism of E. 
Then we have 
$$ 
j(\Delta_a)= \left( 
              \begin{array}{cc}
              1 & \bar{a}\\
              a & \bar{a}a
              \end{array}
              \right) 
$$
Since $\{{\rm id}, \frac{1+\alpha}{2},
 F\frac{1+\alpha}{2},\frac{(a+F)\alpha}{q}\}$ is a basis of ${\cal O} = {\rm End}(E)$,
we see that
$$
   E_{1}, E_{2}, \Delta=\Delta_{\rm id}, \Delta_{\frac{1+\alpha}{2}},
 \Delta_{F\frac{1+\alpha}{2}}, \Delta_{\frac{(a+F)\alpha}{q}}
$$
is a basis of  ${\rm NS}(A)$.
Since $\alpha^{2}=-q $, we see that 
$\alpha$ acts on ${\rm H}^{0}(E, \Omega^{1}_{E})$ as 
multiplication by $\pm \sqrt{-q}$. We can choose $\alpha$ such that
the action $\alpha$ on ${\rm H}^{0}(E, \Omega_{E}^{1})$ is multiplication
by $\sqrt{-q}$. $F$ acts on
${\rm H}^{0}(E, \Omega^{1}_{E})$ as the zero-map. 
Therefore, $\frac{1+ \alpha}{2}, F\frac{1+ \alpha}{2}$ 
and $\frac{(a+F)\alpha}{q}$ act 
on ${\rm H}^{0}(E, \Omega^{1}_{E})$ respectively as multiplication by 
$$ 
\frac {1 + \sqrt{-q}}{2},~ 0,~ \frac{a \sqrt{-q}}{q }.
$$
Since ${\rm H}^{1}(E, {\cal O}_{E})$ is dual to $ {\rm H}^{0}(E, \Omega^{1}_{E})$, 
by Lemma \ref{3.1} the actions of $\frac{1+ \alpha}{2}, F\frac{1+ \alpha}{2}$ 
and $\frac{(a+F)\alpha}{q}$ on ${\rm H}^{1}(E, {\cal O}_{E})$ are respectively
given as multiplication by
$$
\frac {1 - \sqrt{-q}}{2},~ 0,~ -\frac{a \sqrt{-q}}{q }.
$$
Therefore, on the decomposition $(*)$ of the space 
${\rm H}^{1}(A, \Omega^{1}_{A})$
the endomorphisms
$
{\rm id} \times \frac{1+ \alpha}{2}, ~ {\rm id} \times F \frac{1+ \alpha}{2},~ {\rm id} \times \frac{(a+F)\alpha}{q}  
$ 
of $A$ act respectively as multiplication by
$$
(1,\frac{1+ \sqrt{-q}}{2}, \frac{1 - \sqrt{-q}}{2}, \frac{1+q}{4}),\quad
(1, 0, 0, 0),\quad
(1, \frac{ \sqrt{-q}}{q},  - \frac{ \sqrt{-q}}{q},  \frac{a^{2}}{q})
$$
on each direct summand.

We consider the automorphism $\tau$ of $A$ defined by
$$
\begin{array}{cccc}
   \tau : & A = E_{1}\times E_{2} & \longrightarrow & A = E_{1} \times E_{2} \\
             &     (P_{1}, P_{2}) & \mapsto & (P_{2}, P_{1}).  
\end{array}
$$
We denote by $\Omega_{2}$ a basis of 
$ {\rm H}^{1}(E_1, {\cal O}_{E_1})\otimes {\rm H}^{0}(E_2, \Omega^{1}_{E_2})$.
We set $\Omega_{3} = \tau^{*}\Omega_{2}$. Then, $\Omega_{3}$ is a basis of 
${\rm H}^{0}(E_1, \Omega^{1}_{E_1}) \otimes {\rm H}^{1}(E_2, {\cal O}_{E_2})$,
and there exist coefficients $\alpha_{i} \in k$ $(i = 1, 2, 3, 4)$ such that
$$
c_{1}(\Delta)=c_{1}(\Delta_{{\rm id}})= \alpha_{1} \Omega_{1} + \alpha_{2} \Omega_{2} + \alpha_{3} \Omega_{3} + \alpha_{4} \Omega_{4} 
$$
We consider inclusions
$$
\begin{array}{lcccclccc}
           \epsilon_1 : &  E_{1}  & \longrightarrow  &  E_{1}\times E_{2}= A& 
& \epsilon_2 : &  E_{2}  & \longrightarrow  &  E_{1}\times E_{2}= A\\
                &  P     & \mapsto   &     (P, O_{E_2})&  &  &  P     & \mapsto   
&     (O_{E_1}, P)
\end{array}
$$
Then, we have the following diagram induced by $\epsilon_i$.
$$
\begin{array}{ccc}
{\rm Pic}(E_{i})/p{\rm Pic}(E_{i}) & \stackrel{c_{1}}{\longrightarrow} & 
{\rm H}^{1}(E_{i}, \Omega_{E_{i}}^{1}) \\
\uparrow  &       &\uparrow \\
{\rm NS}(A)/p{\rm NS}(A) & \stackrel{c_{1}}{\longrightarrow} & 
{\rm H}^{1}(A, \Omega_{A}^{1}).
\end{array}
$$
Using this diagram, by $\Delta\cdot E_{1}=1$ and $\Delta\cdot E_{2} =1$ we see  
$\alpha_{1}=\alpha_{4}=1$. Since $\tau^{*}\Delta = \Delta$, we also have
$\alpha_{2}= \alpha_{3}$, which we denote by $\alpha$. 

We show now $\alpha \neq 0$.
We consider the natural inclusion $\phi: \Delta \hookrightarrow A = E_{1}\times E_{2}$
and the diagram
$$
\begin{array}{ccc}
{\rm Pic}(\Delta)/p{\rm Pic}(\Delta) & \stackrel{c_{1}}{\longrightarrow} & 
{\rm H}^{1}(\Delta, \Omega_{\Delta}^{1}) \\
\uparrow  &       &\uparrow \\
{\rm NS}(A)/p{\rm NS}(A) & \stackrel{c_{1}}{\longrightarrow} & 
{\rm H}^{1}(A, \Omega_{A}^{1}).
\end{array}
$$
Since $\Delta^2 = 0$, we have $\phi^{*}c_{1}(\Delta) = 0$.
On the other  hand, since $\Delta\cdot E_{1} = \Delta\cdot E_{2} = 1$,
we have $\phi^{*}c_{1}(E_{1})= \phi^{*}c_{1}(E_{2}) \neq 0$.
Therefore, we see $\alpha \neq 0$. Replacing $\Omega_{2}$ by $\alpha\Omega_{2}$,
we may assume $\alpha = 1$.

Summarizing these results, we have
$$
\begin{array}{l}
c_{1}(E_{1}) =\Omega_4, ~c_{1}(E_{2}) =\Omega_1, \\
c_{1}(\Delta)=\Omega_{1} + \Omega_{2} + \Omega_{3} + \Omega_{4}, \\
c_{1}(\Delta_{\frac{1+ \alpha}{2}})= \Omega_{1}+ 
\frac{1 + \sqrt{-q}}{2} \Omega_{2} + \frac{1 - \sqrt{-q}}{2} \Omega_{3} + \frac{1+q}{4} \Omega_{4},\\
c_{1}(\Delta_{F\frac{1+ \alpha}{2}})=\Omega_{1},\\
c_{1}(\Delta_{\frac{(a+F)\alpha}{q}}) =\Omega_{1}+ \frac{a\sqrt{-q}}{q} \Omega_{2} - \frac{a \sqrt{-q}}{q}\Omega_{3} + \frac{a^{2}}q \Omega_{4}.
\end{array}
$$
Since  $2q$ is prime to $p$, there exists an integer 
$\ell$ such that $\ell \equiv \frac{a}{2q}~({\rm mod}~p)$. 
Keeping these notations, we have the following theorem.
\begin{th}\label{mainth}
The kernel $\Ker ~c_{1}$ is $2$-dimensional over ${\bf F}_{p}$, and 
a basis of $\Ker ~c_{1}$ is given by
divisors 
$$
\Delta_{F \frac{1+ \alpha}{2}} - E_{2},~\Delta_{\frac{2+F\alpha}{q}} - \ell \Delta_{\frac{1+ \alpha}{2}} + 2\ell\Delta -
(\ell +1)E_{2} - (1 -q+ 2a)\ell E_{1}.
$$
\end{th}
\proof{With respect to the basis $\langle \Omega_{1}, \Omega_{2}, \Omega_{3}, \Omega_{4} \rangle$, the Chern classes $c_{1}(E_{1})$, $c_{1}(E_{2})$,
$c_{1}(\Delta)$, $c_{1}(\Delta_{\frac{1+ \alpha}{2}})$, $c_{1}(\Delta_{F\frac{1+ \alpha}{2}})$, $c_{1}(\Delta_{\frac{(a+F)\alpha}{q}})$ are respectively represented
as the following vectors:
$$
\begin{array}{l}
{\bf u}_{1} = 
\left(
\begin{array}{c}
0 \\ 0 \\ 0 \\ 1
\end{array}
\right),\quad
{\bf u}_{2} = 
\left(
\begin{array}{c}
1 \\ 0 \\ 0 \\ 0
\end{array}
\right),\quad
{\bf u}_{3} = 
\left(
\begin{array}{c}
1 \\ 1\\ 1 \\ 1 
\end{array}
\right),\\
{\bf u}_{4} = 
\left(
\begin{array}{c}
1 \\ \frac{1 + \sqrt{-q}}{2}  \\ \frac{1 - \sqrt{-q}}{2} \\ \frac{1 + q}{4}
\end{array}
\right),\quad
{\bf u}_{5} = 
\left(
\begin{array}{c}
1 \\ 0 \\ 0 \\ 0
\end{array}
\right),\quad
{\bf u}_{6} = 
\left(
\begin{array}{c}
1 \\ \frac{a\sqrt{-q}}{q}  \\ \frac{- a\sqrt{-q}}{q} \\ \frac{a^2}{q}
\end{array}
\right).
\end{array}
$$
Since ${\bf u}_{1}, {\bf u}_{2}, {\bf u}_{3}, {\bf u}_{4}$ are linearly independent
over ${\bf F}_{p}$ and we have
$$
\begin{array}{l}
  {\bf u}_{5} = {\bf u}_{2}, \\
{\bf u}_{6} = \frac{2a}{q}{\bf u}_{4} -  \frac{a}{q}{\bf u}_{3}
+ (\frac{a}{2q} + 1){\bf u}_{2}
 + (\frac{a}{2q} - \frac{a}{2} + \frac{a^2}{q}){\bf u}_{1},
\end{array}
$$
we see $\dim_{{\bf F}_{p}}{\rm Im}~ c_{1} = 4$. 
Since $\dim_{{\bf F}_{p}}{\rm NS}(A)/p{\rm NS}(A) = 6$, 
we have $\dim_{{\bf F}_{p}}\Ker~ c_{1} = 2$. 
Since $
\{ E_{1}, E_{2}, \Delta, \Delta_{\frac{1+\alpha}{2}}, \Delta_{F\frac{1+\alpha}{2}},
\Delta_{\frac{(a+F)\alpha}{2}}\}$ 
is a basis of $NS(A)/pNS(A)$,
the latter part follows from our construction.
}
Using this theorem, we have the following known corollary (cf.
van der Geer and Katsura \cite{GK}, for instance).
\begin{cor}
Let A be a superspecial abelian surface. Then, ${\rm H}^{1}(A, \Omega_{A}^{1})$
is generated by algebraic cycles.
\end{cor}
\proof{This follows from the fact that 
 ${\bf u}_{1}, {\bf u}_{2}, {\bf u}_{3}, {\bf u}_{4}$ are 
linearly independent also over $k$.
}

\section{Example}
We give here one concrete example.
Assume characteristic $p=3$. Then, there exists  
only one supersingular elliptic curve up to isomorphism
and it is given by
$$
E : y^{2} = x^{3} - x
$$
We consider two automorphisms defined by
$$
\begin{array}{l}
\sigma: x \mapsto x+1, y \mapsto y,\\
\tau:x \mapsto -x, y \mapsto \sqrt{-1}y 
\end{array}
$$
We have a morphism defined by
$$
\begin{array}{cccc}
\pi:& E  &\longrightarrow &{\bf P}^{1}\\
  &(x,y)& \longmapsto  &x
\end{array}
$$
By the result of Ibukiyama (\cite{I}), we have 
$$
{\rm End}(E) = Z+Z\tau+Z\iota\circ\tau+Z\tau\circ\iota\circ\sigma.
$$
Here, $\iota$ is the involution of $E$. 

Let $P$ be the point on ${\bf P}^{1}$ given 
by the local equation $x=0$, and ${\tilde P}$ a point on $E$ 
such that $\pi({\tilde P}) = P$.
We consider an affine open covering $\{U_{0}, U_{1} \}$ of ${\bf P}^{1}$ which is
given by 
$$
 U_{0}=\{x \in {\bf P}^{1} \mid  x \ne \infty \}, U_{1}=\{x \in {\bf P}^{1} \mid x \ne 0 \}.
$$
The divisor $P$ is given by the functions
$$
x ~{\mbox{on}} ~U_{0},~ 1  ~{\mbox{on}} ~U_{1}.
$$
Under the notation, we have the following diagram.
$$
\begin{array}{cccccc}
2{\tilde P} &\in {\rm Pic}(E)/3{\rm Pic}(E)& \stackrel{c_{1}}{\hookrightarrow} & {\rm H}^{1}(E, \Omega^{1}_{E})& \cong & k\\
\uparrow & \uparrow &                   &\uparrow        &   &    \\
P &\in {\rm Pic}({\bf P}^{1})/3{\rm Pic}({\bf P}^{1})&\stackrel{c_{1}}{\hookrightarrow} &{\rm H}^{1}({\bf P}^{1}, \Omega^{1}_{{\bf P}^{1}})& \cong & k
\end{array}
$$
In this diagram, we have $c_{1}(P)=\{\frac{dx}{x}\}$,
and $c_{1}({\tilde P})=\{\frac{dx}{2x}\}$.

We set $A=E_{1} \times E_{2}$ with $E_{1}=E_{2}=E$. 
We consider the Chern class map
$$
c_{1}: {\rm NS}(A)/3{\rm NS}(A) \rightarrow {\rm H}^{1}(A, \Omega^{1}_{A})
\cong {\rm H}^{1}(A, {\cal O}_{A}) \otimes {\rm H}^{0}(A, \Omega^{1}_{A})
$$
We also consider the natural inclusion defined by
$$
\begin{array}{ccccc}
\varphi :& E &\rightarrow & E_{1} \times E_{2} &=A  \\
       & P & \mapsto &(P, O_{E_2}) &
\end{array}
$$
We have  a commutative diagram
$$
\begin{array}{ccc}
{\rm NS}(A)/3{\rm NS}(A) & \stackrel{\varphi^{*}}{\rightarrow}&  {\rm NS}(E)/3{\rm NS}(E)\\
     \downarrow c_{1} &    & \downarrow c_{1} \\
{\rm H}^{1}(A, \Omega^{1}_{A})&  \stackrel{\varphi^{*}}{\rightarrow} & {\rm H}^{1}(E, \Omega_{E}^{1})
\end{array}
$$
Then, we have 
$$
\varphi^{*}(c_{1}(\Delta)) = c_{1}(\varphi^{*}(\triangle))= c_{1}(O_{E})=\{\frac{dx}{2x}\} \neq 0.
$$

We determine the action of ${\rm End} (E)$ on ${\rm H}^{0}(E, \Omega_{E}^{1})$.
A basis of ${\rm H}^{0}(E, \Omega_{E}^{1})$ is given by $\frac{dx}{y}$
and we have
$$
(\iota\circ\sigma)^{*}\frac{dx}{y} = -\frac{dx}{y},\quad
   \tau^{*}\frac{dx}{y} =- \sqrt{-1}\frac{dx}{y}, \quad (\tau\circ\iota\circ\sigma)^{*}\frac{dx}{y} =\sqrt{-1}\frac{dx}{y}.
$$
Since ${\rm H}^{1}(E, {\cal O}_{E})$ is dual to ${\rm H}^{0}(E, \Omega_{E}^{1})$,
the actions of $\iota\circ\sigma$, $\tau$ and $\tau\circ\iota$ are respectively
given as multiplication by
$$
    -1,~\sqrt{-1}, ~-\sqrt{-1}
$$
by Lemma \ref{3.1}.
Since we have
$$
\begin{array}{cl}
{\rm H}^{1}(A, \Omega^{1}_{A}) & \cong {\rm H}^{1}(A, {\cal O}_{A})\otimes 
{\rm H}^{0}(A, \Omega^{1}_{A})\\
      &\cong ({\rm H}^{1}(E_{1}, {\cal O}_{E_1}) \otimes {\rm H}^{0}(E_{1}, {\cal O}_{E_1})) \oplus ({\rm H}^{1}(E_{1}, {\cal O}_{E_1}) \otimes {\rm H}^{0}(E_{2}, \Omega^{1}_{E_2}))  \\
  & \quad \quad \oplus  ({\rm H}^{1}(E_2, {\cal O}_{E_2})  \otimes  {\rm H}^{0}(E_1, \Omega^{1}_{E_1}))  \oplus ({\rm H}^{1}(E_{2}, {\cal O}_{E_2}) \otimes {\rm H}^{0}(E_{2}, \Omega^{1}_{E_2})),
\end{array}
$$
the actions ${\rm id} \times \iota\circ\sigma$, ${\rm id}\times \tau$ and ${\rm id}\times\tau\circ\iota\circ\sigma$ are respectively given as multiplication 
on each summand by
$$
\begin{array}{l}
(1, -1, -1,  1)\\
(1, \sqrt{-1}, -\sqrt{-1}, 1)\\
(1, -\sqrt{-1}, \sqrt{-1}, 1).
\end{array}
$$
By our general theory,
$$
   E_{1}, E_{2}, \Delta, \Delta_{\iota\circ\sigma},
 \Delta_{\tau}, \Delta_{\tau\circ\iota\circ\sigma}
$$
gives a basis of ${\rm NS}(A)$ over ${\bf Z}$. Therefore,
$\Delta+\Delta_{\iota\circ\sigma} + E_{1} + E_{2}$ and 
$ \Delta_{\tau}+ \Delta_{\tau\circ\iota\circ\sigma} + E_{1}  + E_{2}$
are linearly independent divisors in ${\rm NS(A)}/3{\rm NS}(A)$
over ${\bf F}_{3}$. Moreover, considering the actions of the endomorphisms
${\rm id} \times \iota\circ\sigma$, ${\rm id} \times \tau$ and 
${\rm id} \times\tau\circ\iota\circ\sigma$ on ${\rm H}^{1}(A, \Omega_{A}^{1})$
and the commutative diagram
$$
\begin{array}{ccc}
{\rm NS}(A) & \stackrel{f^{*}}{\rightarrow}&  {\rm NS}(A)\\
     \downarrow c_{1} &    & \downarrow c_{1} \\
{\rm H}^{1}(A, \Omega^{1}_{A})&  \stackrel{f^{*}}{\rightarrow} & {\rm H}^{1}(A, \Omega^{1}_{A})
\end{array}
$$
with $f \in {\rm End}(A)$,
we conclude that the Chern classes of these two
divisors are zero. Therefore, we see that 
$$
\Delta+\Delta_{\iota\circ\sigma} + E_{1}  + E_{2},
 \Delta_{\tau}+ \Delta_{\tau\circ\iota\circ\sigma} + E_{1}  + E_{2}
$$
gives a basis of $\Ker ~c_{1}$ over ${\bf F}_{3}$.

\section{An application to Kummer surfaces}
Let $A$ be an abelian surface defined over $k$, and $\iota$ be  the involution 
$x \mapsto \ominus x$.
We denote by ${\tilde{A}}$ the surface made of 16 blowing-ups at 16 two-torsion
points on $A$.
Then, $\iota$ induces the action $\tilde{\iota}$ on ${\tilde{A}}$ and 
${\rm Km}(A) = {\tilde{A}}/\tilde{\iota} $ is the Kummer surface.
We denote by  $\pi : {\tilde{A}} \rightarrow {\rm Km}(A)$ the projection.
A K3 surface $X$ is called supersingular if the Picard number $\rho(X)$
is equal to the second Betti number $b_{2}(X) = 22$. For a supersingular K3
surface, the discriminant of ${\rm NS}(X)$ is equal to the form $-p^{2\sigma_{0}}$
and $\sigma_{0}$ is called an Artin invariant.
We know $1 \leq \sigma_{0} \leq 10$ (cf. Artin \cite{A}). A supersingular K3 surface
with Artin invariant 1 is said to be superspecial. Such a K3 surface is unique
up to isomorphism and is isomorphic to the Kummer surface
${\rm Km}(A)$ such that $A$ is superspecial (cf. Ogus \cite{O}　and Shioda \cite{S}).
We also know that a supersingular K3 surface with $\sigma_{0} = 2$ is isomorphic to
a Kummer surface ${\rm Km}(A)$ such that $A$ is supersingular and
non-superspecial (cf. Ogus \cite{O}).

We have the following commutative diagram:
$$
\begin{array}{ccc}
{\rm NS}({\rm Km}(A))/p{\rm NS}({\rm Km}(A)) & \stackrel{c_{1}}{\longrightarrow} & 
{\rm H}^{1}({\rm Km}(A), \Omega_{{\rm Km}(A)}^{1}) \\
\downarrow  &       &\downarrow \\
{\rm NS}(\tilde{A})/p{\rm NS}(\tilde{A}) & \stackrel{c_{1}}{\longrightarrow} & 
{\rm H}^{1}(\tilde{A}, \Omega_{\tilde{A}}^{1}).
\end{array}
$$
Since we have $2{\rm NS}({\tilde{A}}) \subset \pi^{*}{\rm NS}({\rm Km}(A))
\subset {\rm NS}({\tilde{A}})$ by Shioda \cite{S} and $p \neq 2$,
we see ${\rm NS}({\rm Km}(A))/p{\rm NS}({\rm Km}(A))\cong {\rm NS}(\tilde{A})/p{\rm NS}(\tilde{A})$. Since $\iota$ acts on ${\rm H}^{1}(A, {\cal O}_{A})$
and ${\rm H}^{0}(A, {\Omega}^{1}_{A})$ as multiplication by $-1$.
Since ${\rm H}^{1}(A, {\Omega}^{1}_{A})\cong {\rm H}^{1}(A, {\cal O}_{A})
\otimes {\rm H}^{0}(A, {\Omega}^{1}_{A})$, we see that $\iota$ acts as identity
on ${\rm H}^{1}(A, {\Omega}^{1}_{A})$. Therefore, ${\tilde{\iota}}$
acts as identity on ${\rm H}^{1}({\tilde{A}}, {\Omega}^{1}_{\tilde{A}})$.
Hence, we have ${\rm H}^{1}({\rm Km}(A), \Omega_{{\rm Km}(A)}^{1})\cong
{\rm H}^{1}({\tilde{A}}, {\Omega}^{1}_{\tilde{A}})$.
Summarizing these results, by Theorems \ref{injectivity} and \ref{mainth} 
we have the following theorem.
\begin{th}
For a Kummer surface ${\rm Km}(A)$, let $c_1$ be the Chern class map 
$$
c_{1} : {\rm NS}({\rm Km}(A))/p{\rm NS}({\rm Km}(A)) \longrightarrow
{\rm H}^{1}({\rm Km}(A), \Omega_{{\rm Km}(A)}^{1}) .
$$
Then, we have the following.
\begin{itemize}
\item[$({\rm i})$]  If ${\rm Km}(A)$ is not superspecial, then $c_{1}$ is injective.
\item[$({\rm ii})$] If ${\rm Km}(A)$ is superspecial, then $\dim_{{\bf F}_{p}} 
\Ker ~ c_{1} = 2$.
\end{itemize}
\end{th}
\begin{remark}For a supersingular K3 surface $X$, it is known that the homomorphism 
$$
       {\rm NS}(X)/p{\rm NS}(X)\otimes_{{\bf F}_{p}} k \longrightarrow
{\rm H}^{1}(X, \Omega_X^{1})
$$
induced by $c_{1}$ is not injective  $($cf. Ogus $\cite{O}$$)$. 
In particular, if ${\rm Km}(A)$ is supersingular, 
$$
    {\rm NS}({\rm Km}(A))/p{\rm NS}({\rm Km}(A))\otimes_{{\bf F}_{p}} k 
 \longrightarrow
{\rm H}^{1}({\rm Km}(A), \Omega_{{\rm Km}(A)}^{1}) .
$$
is not injective even if ${\rm Km}(A)$ is not superspecial.
\end{remark}

\vspace{0.5cm}
\noindent
T.\ Katsura: Faculty of Science and Engineering, Hosei University,
Koganei-shi, Tokyo 184-8584, Japan

\noindent
E-mail address: toshiyuki.katsura.tk@hosei.ac.jp
\end{document}